\documentclass[a4paper,10pt]{article}

\title{The Classification of Higher Order Modular Forms and their Cohomology}
\author{David Sim (University of Nottingham)}

\usepackage{amsfonts}
\usepackage{amsmath}
\usepackage{amssymb}
\usepackage{amsthm}

\newtheorem{Theorem}{Theorem}[section]

\newtheorem{Corollary}[Theorem]{Corollary}

\newtheorem{Lemma}[Theorem]{Lemma}
\newtheorem{Proposition}[Theorem]{Proposition}

\theoremstyle{definition}
\newtheorem{Definition}[Theorem]{Definition}

\def\g{\gamma}
\def\d{\delta}
\def\G{\Gamma}
\def\s{\sigma}

\def\C{\mathbb{C}}

\def\BBH{\mathbb{H}}
\def\BBZ{\mathbb{Z}}
\def\BBN{\mathbb{N}}
\def\BBC{\mathbb{C}}
\def\BBR{\mathbb{R}}

\def\CZ{\mathcal{Z}}
\def\CS{\mathcal{S}}
\def\CM{\mathcal{M}}

\def\a{\mathfrak{a}}
\def\b{\mathfrak{b}}
\def\c{\mathfrak{c}}

\newcommand{\ZM}[1]{\mathbf{Z}_{M}^{(#1)}}
\newcommand{\ZS}[1]{\mathbf{Z}_{S}^{(#1)}}

\newcommand{\Zng}[2]{\mathcal{Z}_{i_{#1},...,i_{#2};\, g}}
\newcommand{\Zngj}[2]{\mathcal{Z}_{j_{#1},...,j_{#2};\, g}}
\newcommand{\Znf}[3]{\mathcal{Z}_{i_{#1},...,i_{#2};\, #3}}

\newcommand{\Zpg}[2]{\mathcal{Z'}_{i_{#1},...,i_{#2};\, g}}
\newcommand{\Zpgj}[2]{\mathcal{Z'}_{j_{#1},...,j_{#2};\, g}}
\newcommand{\Zpf}[3]{\mathcal{Z'}_{i_{#1},...,i_{#2};\, #3}}
\newcommand{\Zpog}[2]{\mathcal{Z}'_{-1,1,j_{#1},...,j_{#2};\,g}}

\newcommand{\ZFrac}[2]{\frac{\overline{\ZS{#1}}\bigoplus\ZM{#1}}{\overline{\ZS{#2}}\bigoplus\ZM{#2}}}

\def\cstuff{c_{-1,1,i_{3},...,i_{t-1};\, g}}

\newcommand{\Shuf}[2]{\mathcal{S}_{#1, #2}}

\newcommand{\cusp}[1]{\mathfrak{a}_{#1}}

\newcommand{\Z}[1]{Z_{(#1)}^{1}(\G, P_{k-2})}
\newcommand{\B}[1]{B_{(#1)}^{1}(\G, P_{k-2})}
\newcommand{\Hom}[1]{H_{(#1)}^{1}(\G, P_{k-2})}
\def\Par{\mathrm{par}}

\def\la{\langle}
\def\ra{\rangle}

\begin{document}

\maketitle
\section{Preliminaries}
\subsection{Introduction}
In this paper, we construct an explicit basis for modular forms of all orders. Higher order modular forms are a natural generalization of the classical concept of automorphic forms, and have been to attracting increasing interest in recent years. They have proved to be relevant to problems related to the distribution of modular symbols (\cite{CDO}), to GL(2) $L-$functions (\cite{DKMO},\cite{F}), to percolation theory (\cite{KZ}) and in \cite{DSr} to the non-commutative modular symbols introduced by Manin (\cite{M}), and yielded results including the proof that modular symbols have a normal distribution (\cite{PR}) and the establishment of higher order Kronecker limit formulas (\cite{JO}).\\\\
The classification of higher order modular forms (via the construction of explicit bases) was begun in \cite{DO} with the weight 2 case, and completed in \cite{DS}. The method involved the construction of generalized Poincare series to yield forms satisfying a specific functional equation. These forms were then used via an iterative construction to build up a complete basis.\\\\
The method of this paper mirrors the construction of \cite{DS} very closely. In section 2 we recall the basic results on the Poincare series from \cite{DS}. Section 3 then uses these series to construct some additional non-cuspidal forms. In section 4 we apply the iterative procedure of \cite{DS} to produce a set of forms, which we then prove forms a basis. In the final section, we generalize the cohomological results of \cite{DO} to all orders.
%
%
%
%
%
%
\subsection{Definitions}
We begin by restating some definitions and results from \cite{DS}.\\\\
Let $\G\subseteq\textrm{PSL}(2,\BBR)$ be a Fuchsian group of the first kind acting on the upper half plane $\BBH$ with compact quotient $\G\backslash\BBH$ of genus $g.$ We assume that there are $m\geq2$ inequivalent cusps. We fix a fundamental domain $\mathfrak{F}$ and representitives $\cusp{1},...,\cusp{m}$ of the inequivalent cusps of $\bar{\mathfrak{F}}.$ As in \cite{I}, we have scaling matrices $\sigma_{\cusp{i}}$ taking neighbourhoods of $i\infty$ to neighbourhoods of $\cusp{i}$. Writing $\G_{\cusp{i}}$ for the stabilizers $\mathrm{stab}_{\G}(\cusp{i}),$ we have
\[ \sigma_{\a_{i}}^{-1}\G_{\a_{i}}\sigma_{\a_{i}} = \G_{\infty} = \{ \pm \Bigl(\begin{array}{c c}1&m\\0&1\end{array}\Bigr)|m\in\BBZ\}.\]
We write $\pi_{\cusp{i}}$ for the generator of $\G_{\a_{i}}$ given by $\sigma_{\a_{i}}\Bigl(\begin{array}{c c}1&m\\0&1\end{array}\Bigr)\sigma_{\a_{i}}^{-1}$\\\\
As usual, the slash operator $|_{k}$ defines an action of $\textrm{PSL}(2,\BBR)$ on functions $f:\BBH\rightarrow\BBC$ by
\[(f|_{k}(\g)(z) = f(\g z)j(\g,z)^{-k}\]
where $j(\g,z) = cz+d$ for $\g = \Bigl(\begin{array}{cc}*&*\\c&d\end{array}\Bigr),$ and extend to $\BBZ[\textrm{PSL}(2, \BBR)]$ by linearity.\\\\
Before stating the definition of higher order modular forms and higher order cusp forms, we will state and label various conditions used in the definition. For a function $f: \BBH \to \C$ we have conditions
\begin{itemize}
\item $f$ is holomorphic on $\BBH$ (``holomorphicity'').\label{HCond}
\item $f|_{k}(\pi - 1) = 0$ for all $\pi$ parabolic in $\G$ (``parabolic invariance''). \label{PCond}
\item $f_{k}(\g - 1) \in R$ for $R$ a set of functions $\C \to \BBH$ (``modularity with periods in $R$'').
\item For each cusp $\a,$ $(f|_{k}\sigma_{\a})(z) \ll e^{-cy}$ as $y\to\infty$ uniformly in $x$ with $c >0$ (``vanishing at the cusps'').
\item For each cusp $\a,$ $(f|_{k}\sigma_{\a})(z) \ll c$ as $y\to\infty$ uniformly in $x$ with $c$ constant (``boundedness at cusps'').
\end{itemize}
\begin{Definition}
We write $S^{t}_{k}(\G)$ for the space of cusp forms of weight $k$ and order $t$ for $\G$. For $t = 0,$ this is the set $\{0\}.$ Otherwise it is the space of functions satisfying holomorphicity, parabolic invariance, modularity with periods in $S_{k}^{t-1}(\G)$ and vanishing at cusps.\\\\
Similarly we write $M^{t}_{k}(\G)$ for the space of modular forms of weight $k$ and order $t$ for $\G$. For $t = 0,$ this is again defined to be $\{0\}.$ Otherwise it is the space of functions satisfying holomorphicity, parabolic invariance, modularity with periods in $M_{k}^{t-1}(\G)$ and boundedness at cusps.
\end{Definition}
We fix once and for all bases $\CS_{2} := \{f_{1},...,f_{g}\}$ for $S_{2}(\G)$ and $\CS_{k}$ for $S_{k}(\G).$
\section{Constructing Cusp Forms from Series}
\subsection{Poincare Series}
We will first state, in this and the next section, some of the key definitions and results of \cite{DS}. For proofs, refer to that paper.
\begin{Definition}Given $i_{1},...,i_{t}\in \{1,...,g\}$ and a cusp $\a,$ set
 \[F_{i_{1},...,i_{t}}^{\a}(z) = f_{i_{1}}(z)\int_{z}^{\a}f_{i_{2}}(t_{1})\int_{t_{1}}^{\a}f_{i_{3}}(t_{2})...\int_{t_{t-2}}^{\a}f_{i_{t}}(t_{t-1})dt_{t-1}...dt_{1}\]
Furthermore, write $A_{t}$ for $\{F_{i_{1},...,i_{t}}:i_{j}\in\{1,...,g\}\}.$
\end{Definition}
\begin{Definition}For $m \geq 0,$ $k\in 2\BBZ,$ $\a$ a cusp and $f\in S_{2}(\G)$, set
 \[Z_{\a m}(z,s,1,k;\overline{f}):=\sum_{\g\in\G_{\a}\\\G}\overline{\Bigl(\int_{\a}^{\g\a}f(w)dw\Bigr)}Im(\s_{\a}^{-1}\g z)^{s}e(m\s_{\a}^{-1}\g z)\epsilon(\s_{\a}^{-1}\g,z)^{-k}.\]
\end{Definition}
\begin{Theorem}
There exists $\d_{\G}>0$ such that for any $f\in A_{t,k}$ $Z_{\a m}(z,s,1,k,\overline{f})$ admits a meromorphic continuation to $Re(s) > 1 - \d_{\G}.$ The only possible pole is at $s = 1$ and it can occur only when $k\le 0.$ For $k=0$ it is simple. For $k\geq 2$ and $m\neq 0,$ $Z_{\a m}(z,s,1,k;\overline{f})\ll y_{\mathfrak{F}}(z)^{1/2}.$ For $k\geq 2,$ $Z_{\a 0 }(z,s,1,k;\overline{f}) \ll y_{\mathfrak{F}}(z)^{\s}.$ For $k=2,$ $Z_{\a 0}(z,1,1,2,\overline{f})\ll y_{\mathfrak{F}}(z)^{1/2}.$ The implied constants are independent of $z$ in all cases.\label{maintheorem}
\end{Theorem}
\subsection{Explicit Construction of Cusp Forms}
\begin{Definition}We set
\[Z_{\a m}(z,s;\overline{f}) := y^{-1}Z_{\a m}(z, s+1, 1, 2; \overline{f}).\]
\end{Definition}
It is easy to see that
\begin{multline}Z_{\a m}(\cdot,0;\overline{F}_{i_{1},...,i_{t-1}})|_{2}(\g-1) = \\
\Bigl(\overline{\int_{\a}^{\g^{-1}\a}F_{i_{1},...,i_{t-1}}}\Bigr)P_{\a m}+\sum_{r=1}^{t-2}\Bigl(\overline{\int_{\a}^{\g^{-1}\a}F_{i_{1},...,i_{r}}}\Bigr)Z_{\a m}(\cdot,0;\overline{F}_{i_{r+1},...,i_{t-1}})\label{Zhit}\end{multline}
and thus that these series will obey appropriate functional equations. Their vanishing at cusps is guaranteed by (\ref{maintheorem}).
\cite{DO} analyzes the nonholomorphic part of $Z_{\a m}(z,s;\overline{F}_{i_{1},...,i_{t-1}})$ in order to produce holomorphic linear combinations of them.  In particular, one may construct, for any $i_{1},...,i_{t-1}\in \{1,...,g\}$ and $f \in \CS_{k}$ not satisfying $f = f_{i_{t-1}} = f_{1},$ a function
$\mathcal{Z}_{-i_{1},...,-i_{t-1};f}(\cdot) \in S^{t}_{k}$ satisfying 
\begin{multline*}
\mathcal{Z}_{-i_{1},...,-i_{t-1};f}|_{2}(\g_{1}-1)...(\g_{t-1}-1)\\
= \left\{\begin{array}{ll}\overline{\la f_{i_{1}},\g_{1}\ra}...\overline{\la f_{i_{t-1}}, \g_{t-1}\ra}f & f \neq f_{i_{t}}\\
\overline{\la f_{i_{1}},\g_{1}\ra}...\overline{\la f_{i_{t-2}}, \g_{t-2}\ra}\Bigl(\overline{\la f_{i_{t-1}}, \g_{t-1}\ra}f-\overline{\la f_{i_{1}}, \g_{1}\ra}f_{1}\Bigr) & f = f_{i_{t}}
\end{array}\right.
\end{multline*}
Furthermore, for any $i_{1},...,i_{t-1}\in\{1,...,g\}$ and $f\in \CS_{k},$ one can construct a function $\mathcal{Z}'_{-i_{1},...,-i_{t-1};f}$
that is holomorphic, parabolically invariant, vanishes at every cusp except in the case $t=2$ where it can have polynomial growth at $\a_{m},$ and which satisfies the functional equation
\[\mathcal{Z}'_{-i_{1},...,-i_{t-1};f}|_{2}(\g_{1}-1)...(\g_{t-1}-1) = \overline{\la f_{i_{1}},\g_{1}\ra}...\overline{\la f_{i_{t-1}}, \g_{t-1}\ra}f.\]
In cases where $f \neq f_{i_{t-1}}$ (including all cases where $k>2$), this function is simply $\mathcal{Z}_{-i_{1},...,-i_{t-1};f}$ - it is only when $f = f_{i_{t}}$ that a new construction is needed.
%
%
%
%
%
%
\section{Extension to the Modular Forms}
\subsection{The Weight 2 Case}
Our first goal in this discussion is to construct similar functions for $f\in \CM_{2}.$ For this, we will need some preliminary results.
\begin{Lemma}For $\a$ and $\b$ cusps
  \begin{multline*}\int_{z}^{\a}F_{i_{1},...,i_{t}}^{\a}(t)dt-\int_{z}^{\b}F_{i_{1},...,i_{t}}^{\b}(t)dt\\ = \int_{\b}^{\a}F_{i_{1},...,i_{t}}^{\a}(t)dt - \sum_{r=1}^{t-1}\int_{z}^{\b}F_{i_{1},...,i_{r}}^{\b}(t)dt\cdot\int_{\b}^{\a}F_{i_{r+1},...,i_{t}}^{\a}(t)dt\end{multline*}
 \end{Lemma}
\begin{proof}
 Induction.
\end{proof}
We now define the functions that we will use to study the residues of the forms we are interested in.
\begin{Definition}Given a cusp $\a$, set $S^{\a}_{i_{r}} = \overline{\int_{z}^{\a}f_{i_{r}}(t)dt}$ and $S^{\a}_{i_{r+1},i_{r}} = 1.$ Then define recursively
\[S^{\a}_{i_{1},...,i_{t}} = \sum_{r=1}^{t}\overline{\int_{z}^{\a}F_{i_{1},...,i_{r}}^{\a}(t)dt}\cdot S_{i_{r+1},...,i_{t}}^{a}(z).\]
\end{Definition}
\begin{Proposition}For $\a$ and $\b$ cusps,
\[
S^{\a\b}_{i_{1},...,i_{t}}(z) := S^{\a}_{i_{1},...,i_{t}}(z) - S^{\b}_{i_{1},...,i_{t}}(z)= \sum_{r=1}^{t}\overline{\int_{\b}^{\a}F_{i_{1},...,i_{r}}^{\a}(t)dt}\cdot S^{\a}_{i_{r+1},...,i_{t}}(z)
\] for all $t.$
\end{Proposition}
\begin{proof}
The case $t = 1$ is obvious. Now work inductively. By rearranging we see that
\begin{multline*}S^{\a\b}_{i_{1},...,i_{t}} = \sum_{r=1}^{t}\Bigl(\overline{\int_{z}^{\b}F^{\b}_{i_{1},...,i_{r}}(t)dt}\cdot S^{\a\b}_{i_{r+1},...,i_{t}}\\
+\Bigl(\overline{\int_{z}^{\a}F^{\a}_{i_{1},...,i_{r}}(t)dt} - \overline{\int_{z}^{\b}F^{\b}_{i_{1},...,i_{r}}(t)dt}\Bigr)S_{i_{r+1},...,i_{t}}^{\a}\Bigr)\end{multline*}
We can now apply the inductive hypothesis and the previous lemma to this expression to arrive at
\begin{multline*}
S^{\a\b}_{i_{1},...,i_{t}} = \sum_{r=1}^{t}\sum_{s=r+1}^{t}\overline{\int_{z}^{\b}F^{\b}_{i_{1},...,i_{r}}(t)dt} \overline{\int_{\b}^{\a}F_{i_{r+1},...,i_{s}}^{\a}(t)dt}\cdot S^{\a}_{i_{s+1},...,i_{t}}\\
 -\sum_{r=1}^{t}\sum_{s=1}^{r-1}\overline{\int_{z}^{\b}F_{i_{1},...,i_{s}}^{\b}(t)dt}\overline{\int_{\b}^{\a}F_{i_{s+1},...,i_{r}}^{\a}(t)dt}\cdot S_{i_{r+1},...,i_{t}}^{\a}\\
+\sum_{r=1}^{t}\overline{\int_{\b}^{\a}F_{i_{1},...,i_{r}}^{\a}(t)dt}\cdot S_{i_{r+1},...,i_{t}}^{\a}.
\end{multline*}
Cancelling terms in the first two sums we complete the inductive step, thereby establishing the propostion.
\end{proof}
We can see from the series expression that
\[\frac{d}{d\overline{z}}Z_{\a 0}(z,0;\overline{F}_{i_{1},...,i_{t-1}}) = \frac{i}{2y^{2}}\textrm{Res}_{s=1}Z_{\a 0}(z,s,1,0; \overline{F}_{i_{1},...,i_{t-1}})\] and also by (\cite{DS}, Lemma 3.8) that
\[\textrm{Res}_{s=1}Z_{\a 0}(\cdot,s,1,0;\overline{F}_{i_{1},...,i_{t}}) = \frac{1}{V}S^{a}_{i_{1},...,i_{t}}.\]
Taking our set $\{\a_{1},...,\a_{m}\}$ of inequivalent cusps, we now define functions $f_{\a} := P_{\a 0} - P_{\a_{m}0}$ for $\a \neq \cusp{m}$; this gives a basis for the space of Eisenstein series in $M_{2}(\G).$ We will write $\CM_{2}$ for the basis for $M_{2}(\G)$ consisting of this basis together with $\CS_{2}.$\\\\
For $i_{1},...,i_{t-1} \in [1,g]$ and $\a \neq \a_{m}$ we then set \begin{multline}\CZ_{-i_{1},...,-i_{t-1}; f_{\a}}(z) = Z_{\a 0}(z,0;\overline{F}_{i_{1},...,i_{t-1}}^{\a_{i_{t}}}) - Z_{\a_{m} 0}(z,0;\overline{F}_{i_{1},...,i_{t-1}}^{\a_{m}})\\ - \sum_{r=1}^{t}\overline{\int_{\a_{m}}^{\a_{i_{t}}}F_{i_{1},...,i_{r}}^{\a_{i_{t}}}(t)dt}\cdot Z_{\a_{i_{t}} 0}(z,0;\overline{F}^{\a_{i_{t}}}_{i_{r+1},...,i_{t}}(z)).\label{Zpseries}\end{multline}
We also set $\CZ'_{-i_{1},...,-i_{t-1};f_{\a}} = \CZ_{-i_{1},...,-i_{t-1};f_{\a}}.$
\begin{Proposition}$\CZ_{-i_{1},...,-i_{t-1};f_{\a}}\in M^{t}_{2}(\G).$\label{allneg}\end{Proposition}
\begin{proof}
We can apply the previous result to see that $\frac{d}{d\overline{z}}\CZ_{-i_{1},...,-i_{t-1};f_{\a}} = 0.$ We can also see from the functional equation for $Z_{\a m}(z,0;\overline{F}_{i_{1},...,i_{t-1}})$ that
$\CZ_{-i_{1},...,-i_{t-1};f_{\a}}$ satisfies the functional equation and parabolic invariance conditions.\\\\
By (\ref{maintheorem}), $Z_{\a 0}(z,s,1,k;\overline{f}) \ll y_{\mathfrak{F}}(z)^{\sigma}$ for $k\geq 2$, and so $y^{-1}Z_{\a 0}(z,1,1,2;\overline{f}) \ll y_{\mathcal{F}}(z)^{0}.$ By (\ref{Zhit}), this means that both $Z_{\a 0}(\cdot,0;\overline{F}_{i_{1},...,i_{t-1}}^{\a})$ and $Z_{\a 0}(\cdot,0;\overline{F}_{i_{1},...,i_{t-1}}^{\a})|_{2}(\g-1)$ satisfy this growth condition, and thus that $Z_{\a 0}(\cdot,0;\overline{F}_{i_{1},...,i_{t-1}}^{\a})|_{2}\g$ does as well.
\end{proof}
%
%
%
%
\subsection{Higher Weights}
For $k > 2$, the construction is simpler. Observe first that a basis for the Eisenstein series in $M_{k}(\G)$ is given by the functions $f_{\a}:= P_{\a 0}(z)_{k}$ for $\a \in \{\a_{1},...,\a_{m}\}.$ Adding these functions to $\CS_{k}$ gives us a basis $\CM_{k}$ for $M_{k}(\G).$\\\\
We consider
\begin{multline*}y^{-k/2}Z_{\a 0}(z,s+k/2,1,k;\overline{f})\\
= \sum_{\g\in\G_{\a}\backslash\G}\Bigl(\overline{\int_{\a}^{\a}f(w)dw}\Bigr)\Im (\sigma_{\a}^{-1}\g z)^{s}j(\sigma_{\a}^{-1}\g,z)^{-k}.\end{multline*}
By (\ref{maintheorem}) this extends to an analytic function of $s$ for $\Re(s) > 1 - k/2 -\d_{\G},$ and it is easy to see that
\begin{multline*}y^{-k/2}Z_{\a 0}(z,k/2,1,k;\overline{F}_{i_{1},...,i_{t-1}})|_{k}(\g-1)=\\
\Bigl(\overline{\int_{\a}^{\g^{-1}\a}F_{i_{1},...,i_{t-1}}}\Bigr)P_{\a0}(z)_{k} + \sum_{r=1}^{t-2}\Bigl(\overline{\int_{\a}^{\g^{-1}\a}F_{i_{1},...,i_{r}}}\Bigr)y^{-k/2}Z_{\a 0}(z,k/2,1,k;\overline{F}_{i_{r+1},...,i_{t-1}}).\end{multline*}
Now, for $\Re(s)$ large, we can differentiate term by term to get
\[ \frac{d}{d\overline{z}} y^{-k/2}Z_{\a 0}(z,s+k/2,1,k;\overline{F}_{i_{1},...,i_{t-1}}) = \frac{is}{2y^{1+k/2}}Z_{\a 0}(z,s+k/2,1,k-2;\overline{F}_{i_{1},...,i_{t-1}})\]
By \ref{maintheorem}, $Z_{\a 0}(z,s+k/2,1,k-2;\overline{F}_{i_{1},...,i_{t-1}})$ is holomorphic at $s=0$ if $k>2,$ and so by comparing analytic continuations, we find that $ y^{-k/2}Z_{\a 0}(z,k/2,1,k;\overline{F}_{i_{1},...,i_{t-1}})$ is holomorphic in $z$.\\\\
Now, for any $\a$ we can set
\[\CZ_{-i_{1},...,-i_{t-1};f_{\a}} = Z_{\a 0}(z,k/2,1,k;\overline{F}_{i_{1},...,i_{t-1}})\]
for $i_{1},...,i_{t-1}\in\{1,...,g\}$ and again, we define $\CZ'$ identically:
\[\CZ'_{-i_{1},...,-i_{t-1};f_{\a}} = \CZ_{-i_{1},...,-i_{t-1};f_{\a}}.\]

%
%
%
%
%
%
\section{Constructing the Basis}
\subsection{Combinatorial Preliminaries}
Given the functions defined in the last section, the construction of the basis will now be a primarily combinatorial exercise.
\begin{Definition}For $r, t \in \BBN$ with $r < t,$ a shuffle of type $(r,t)$ is a pair $(\phi, \psi)$ of order preserving maps 
\[\phi:\{1,...,r\}\to \{1,...,t-1\}\]
\[\psi:\{r+1,...,t-1\}\to\{1,...,t-1\}\]
whose images are disjoint and complementary. For convenience, we will denote by $(\phi,\psi)_{0}$ the shuffle such that $\phi(i) = i$ for $1 \le i \le r$ and $\psi(i) = i$ for $r+1\le i \le t-1$, where $r$ and $t$ should be obvious from the context. We write $\mathcal{S}_{r,t}$ for the set of shuffles of type $(r,t).$ The following result is proved in \cite{CD}.
\end{Definition}
\begin{Proposition}Given $F,G: \BBH \to \BBC$ satisfying $F|_{k_{1}}(\g_{1}-1)...(\g_{s}-1) = 0$ and $G|_{k_{2}}(\d_{1}-1)...(\d_{t}-1) = 0$ for all $\g_{1},...,\g_{s},\d_{1},...,\d_{t}\in\G,$
\begin{multline*}F\cdot G|_{k_{1}+k_{2}}(\g_{1}-1)...(\g_{s+t-2}-1)\\
= \sum_{(\phi,\psi)\in \Shuf{s-1}{s+t-1}} F|_{k_{1}}(\g_{\phi(1)}-1)...(\g_{\phi(s-1)}-1)\cdot G|_{k_{2}}(\g_{\psi(s)}-1)...(\g_{\psi(s+t-1)}-1)\end{multline*}
\label{shuffles}\end{Proposition}
We now define index sets with which to label the basis elements we create. We set
\begin{eqnarray*}J_{t,k} &= &\bigl\{(i_{1},...,i_{t-1};g):i_{1},...,i_{t-1}\in\{\pm1,...,\pm g\}, g\in \CM_{k}\bigr\}\\
I_{t,k}& = &\bigl\{(i_{1},...,i_{t-1};g)\in J_{t}: \nexists\: j \textrm{ with } -i_{j} = i_{j+1} = 1\textrm{ and we do not have } g = f_{-i_{t-1}} = f_{1}\bigr\}.\end{eqnarray*}
Note that the condition on $g$ and $f_{i_{t-1}}$ is only relevant in the case $k=2.$ We also define $\mathcal{A}_{t,k}$ to be the linear span over $(i_{1},...,i_{t-1};g)\in J_{t,k} - I_{t,k}$ of the maps
\[\phi:(\g_{1},...,\g_{t-1})\mapsto \la f_{1_{1}},\g_{1}\ra...\la f_{i_{t-1}},\g_{t-1}\ra g.\]
\subsection{The Iterative Construction}
%
%
%
%
%
%
%
%
\begin{Theorem}
For any $(i_{1},...,i_{t-1};g) \in I_{t,k}$ there exists a function $\CZ_{i_{1},...,i_{t-1};g}\in M_{k}^{t}(\G)$ satisfying the functional equation
\[\CZ_{i_{1},...,i_{t-1};g}|_{k}(\g_{1}-1)...(\g_{t-1}-1) = \langle f_{i_{1}}, \g_{1} \rangle ... \langle f_{i_{t-1}}, \g_{t-1} \rangle g + \phi(\g_{1},...,\g_{t-1})\label{PhiFE}\] for some $\phi \in A_{t,k}.$\label{mainthe}
\end{Theorem}
\begin{proof}
First, in the case $t=1$ we set $\CZ_{g} = g.$\\\\
We then can proceed with the iterative construction precisely as in the cuspidal case. Specifically, for any $s<t$ and any $k$ we assume the existence of $\CZ_{i_{1},...,i_{s-1};g}$ for all $(i_{1},...,i_{s-1};g)\in I_{s,k}.$ We then proceed iteratively for each $k$, constructing $\Zng{1}{t-1}$ for $(i_{1},...,i_{t-1};g)\in I_{t,k}$ that satisfy $i_{r}>0$ and $i_{r+1},...,i_{t-1} < 0,$ first for $r = t-1$ and then inductively for lower $r.$\\\\
For $(i_{1},...,i_{t-1};g)$ such that $i_{t-1} > 1,$ we set
\[\Zng{1}{t-1}(z)=g(z)\int_{i}^{z}\Znf{1}{t-2}{f_{i_{t-1}}}(w)dw.\]
By \ref{shuffles} and the inductive hypothesis, this satisfies
\[\Zng{1}{t-1}|_{k}(\g_{1}-1)...(\g_{t-1}-1) = \la f_{i_{1}},\g_{1}\ra...\la f_{i_{t-1}},\g_{t-1}\ra g +\phi(\g_{1},...,\g_{t-1})\]
as desired.\\\\
Now assume that for each $q>r$ we have constructed $\Zngj{1}{t-1}$ for all $(j_{1},...,j_{t};g)\in I_{t,k}$ with $j_{q} > 0,$ $j_{q+1},...,j_{t-1} < 0.$ Now, given $(i_{1},...,i_{t})\in I_{t,k}$ with $i_{r}>0$ and $i_{r+1},...,i_{t-1}<0,$ \ref{shuffles} implies that
\begin{multline}
\Bigl(\Zng{r+1}{t-1}\int_{i}^{z}\Znf{1}{r-1}{f_{i_{r}}}(w)dw\Bigr)|_{k}(\g_{1}-1)...(\g_{t-1}-1)=\\
\sum_{(\phi,\psi)\in\Shuf{r}{t}} \la f_{i_{1}},\g_{\phi(1)}\ra...\la f_{i_{r}},\g_{\phi(r)}\ra\la f_{i_{r+1}},\g_{\psi(r+1)}\ra...\la f_{i_{t-1}},\g_{\psi(t-1)}\ra g\\ + \phi(\g_{1},...,\g_{t-1})\label{inthit}\end{multline}
for some $\phi \in A_{t,k}.$ Consider an individual term in the sum corresponding to $(\phi,\psi)\neq (\phi,\psi)_{0}$  and take $(j_{1},...,j_{t-1})$ such that
\begin{multline*}\la f_{i_{1}},\g_{\phi(1)}\ra...\la f_{i_{r}},\g_{\phi(r)}\ra\la f_{i_{r+1}},\g_{\psi(r+1)}\ra...\la f_{i_{t-1}},\g_{\psi(t-1)}\ra =\\ \la f_{j_{1}},\g_{1} \ra...\la f_{j_{t-1}},\g_{t-1}\ra.\end{multline*}
Because $(\phi,\psi)\neq (\phi,\psi)_{0},$ we must have $\phi(r) \geq r+1,$ and thus since $i_{r} > 0,$ there must be a $q>r$ such that $j_{q}>0$ and $j_{q+1},...,j_{t-1}<0.$ If $(j_{1},...,j_{t-1};g)\in I_{t,k},$ then the inductive hypothesis means that there is a $\Zngj{1}{t-1}$ already constructed such that $\Zngj{1}{t-1}|_{k}(\g_{1}-1)...(\g_{t-1}-1) = \la f_{i_{1}},\g_{\phi(1)}\ra...\la f_{i_{t-1}},\g_{\psi(t-1)}\ra g.$ Otherwise, $\la f_{i_{1}},\g_{\phi(1)}\ra...\la f_{i_{t-1}},\g_{\psi(t-1)}\ra g \in A_{t,k}.$\\\\
Thus if we write $\mathcal{B}(i_{1},...,i_{t-1};g)$ for the sum of all the $\Zngj{1}{t-1}$ with $(j_{1},...,j_{t-1};g)\in I_{t,k}$ corresponding to $(\phi,\psi)\neq (\phi, \psi)_{0},$
then we can define
\[\Zng{1}{t-1} = \Zng{r+1}{t-1}\int_{i}^{z}\Znf{1}{r-1}{f_{i_{r}}}(w)dw - \mathcal{B}(i_{1},...,i_{t-1};g).\]
The only term in  $\Bigl(\Zng{r+1}{t-1}\int_{i}^{z}\Znf{1}{r-1}{f_{i_{r}}}(w)dw\Bigr)|_{k}(\g_{1}-1)...(\g_{t-1}-1)$ that is not either equal to $\Zngj{1}{t-1}|_{k}(\g_{1}-1)...(\g_{t-1}-1)$ for one of the terms in $\mathcal{B}(i_{1},...,i_{t-1};g)$ or in $A_{t,k}$ is the term corresponding to $(\phi, \psi)_{0},$ and so
\[\Zng{1}{t-1}|_{k}(\g_{1}-1)...(\g_{t-1}-1) = \la f_{i_{1}},\g_{1}\ra...\la f_{i_{t-1}},\g_{t-1}\ra g +\phi(\g_{1},...,\g_{t-1})\] for some $\phi \in A_{t,k},$ as desired.\\\\
%
%
The case where all $i_{1},...,i_{t-1} < 0$ cannot be constructed from an integral in the above manner. Here, we use the forms $\Zng{1}{t-1}$ constructed earlier.\\\\
These functions obey the required functional equation by construction. Since each $\Zng{1}{r-1}$ is a weight k order $r$ cusp form, their integrals are weight 0 order $r+1$ cusp forms, and so the products are all weight k order $t$ modular forms. 
\end{proof}
%
%
%
%
%
%
%
\begin{Theorem}
For any $(i_{1},...,i_{t-1};\,g) \in J_{t}$ there exists a function $\Zpg{1}{t-1}$ satisfying the holomorphicity condition, boundedness at cusps, and the functional equation
\[\Zpg{1}{t-2}|_{k}(\g_{1}-1)...(\g_{t-1}-1) = \langle f_{i_{1}}, \g_{1} \rangle ... \langle f_{i_{t-1}}, \g_{t-1} \rangle g.\label{NoPhiFE}\] Furthermore, the $\Zpg{1}{t-1}$ are invariant under $\pi_{\a_{i}}$ for $i < m,$ and if $i_{t} \le g,$ vanish at all cusps except for $\a_{m}.$ 
\end{Theorem}
\begin{proof}
The construction here is almost identical to the construction of $\Zng{1}{t-1},$ so the details will not be given. We start as before with $\CZ'_{g} = g.$ The iterative construction proceeds along the same lines - we define
\[\Zpg{1}{t-1} = \Zpg{r+1}{t-1}\int_{i}^{z}\Zpf{1}{r-1}{f_{i_{r}}}(w)dw - \mathcal{C}(i_{1},...,i_{t-1};g)\]
where $\mathcal{C}(i_{1},...,i_{t-1};g)$ is now a sum over $(j_{1},...,j_{t-1};g) \in J_{t,k}$ of terms $\Zpgj{1}{t-1}$ corresponding to the shuffles in the analogue of (\ref{inthit}). This larger index set gives us a term for every shuffle, and consequently no $\phi$ is needed in the functional equations.\\\\
The main difference comes from the definition of $\mathcal{Z}'_{i_{1},...,i_{t-1};\,g}$ for $i_{1},...,i_{t-1}<0.$ Here, if $g\in S_{k}(\G),$ we use the $\mathcal{Z}'_{i_{1},...,i_{t-1};\,g}$ constructed in \cite{DS}. If $g \in M_{k}(\G)-S_{k}(\G),$ we use the forms constructed earlier in this paper.\\\\
The desired properties of these series all follow from the construction or from previous theorems.
\end{proof}
%
%
%
%
%
%
\subsection{Some Technical Results}
\begin{Lemma}
For all $t \geq 3$ and for any $(i_{1},...,i_{t-1},g)\in J_{t}$ with $i_{1} < 0,$
\begin{multline*} \Zpg{1}{t-1}|_{k}(\pi_{\cusp{m}}-1)(\g_{3}-1)...(\g_{t-1}-1) \\= \left\{\begin{array}{cc}\int_{i}^{\pi_{\cusp{m}}i}\CZ'_{-1;\,f_{1}}(w)dw\la f_{i_{3}},\g_{3}\ra\la f_{i_{t-1}},\g_{t-1}\ra g & \textrm{ if }(i_{1},i_{2}) = (-1,1)\\
0& \textrm{otherwise}\end{array}\right.\end{multline*}\label{Zpcusp}
\end{Lemma}
\begin{proof}
We proceed inductively, following the order of the iteration by which the functions were constructed. The base case, $t = 3$, is a simple calculation. We should also note that in the case $t=2,$ it follows directly from the definitions that $\CZ'_{i_{1};\,g}|_{2}(\pi_{\a_{m}} - 1) = 0$ for all $i_{1}$ and $g.$\\\\
Assuming the result for $s < t,$ if $(i_{1},...,i_{t-1})$ has $i_{t-1}> 0,$ then $\Zpg{1}{t-1} = g\int_{i}^{z}\Zpf{1}{t-2}{f_{i_{t-1}}}(w)dw$ by definition and the result follows by the inductive hypothesis. Thus for induction we will assume further that for $r'>r$ the result holds for any $(i_{1},...,i_{t-1})$ with $i_{r'}>0$ and $i_{r'+1},...,i_{t-1}<0.$\\\\
We now consider $(i_{1},...,i_{t-1};\,g)\in J_{t,k}$ with $i_{1}<0,$ $i_{r}>0$ and $i_{r+1},...,i_{t-1}<0.$ Recall that \[\Zpg{1}{t-1} = \Zpg{r+1}{t-1}(z)\int_{i}^{z}\Zpf{1}{r-1}{f_{i_{r}}}(t)dt - \sum\Zpgj{1}{t-1}\]
where the sum is over all $(j_{1},...,j_{t-1})$ such that \begin{multline*}\Zpgj{1}{t-1}|_{k}(\g_{1}-1)...(\g_{t-1}-1) = \la f_{j_{1}},\g_{1}\ra...\la f_{j_{t-1}},\g_{t-1}\ra g \\ =\la f_{i_{1}},\g_{\phi(1)}\ra...\la f_{i_{r}},\g_{\phi(r)}\ra \la f_{i_{r+1}},\g_{\psi(r+1)}\ra...\la f_{i_{t-1}},\g_{\psi(t-1)}\ra g \end{multline*} for some shuffle $(\phi, \psi) \in \Shuf{r}{t}$ with $(\phi,\psi)\neq (\phi,\psi)_{0}$.\\\\
We note first that since $i_{r+1},...,i_{t-1} < 0,$ the inductive hypothesis and (\ref{shuffles}) imply that \begin{multline*}\Bigl(\Zpg{r+1}{t}(z)\int_{i}^{z}\Zpf{1}{r-1}{f_{i_{r}}}(t)dt\Bigr)|_{k}(\pi_{\cusp{m}}-1)(\g_{3}-1)...(\g_{t-1}-1) = \\
\displaystyle{\sum_{\substack{(\phi,\psi)\in\Shuf{r}{t}\\\phi(1)=1,\phi(2) = 2}}}\int_{i}^{\pi_{\cusp{m}}i}\CZ'{-1;\,f_{1}}(t)dt\la f_{i_{3}},\g_{\phi(3)}\ra...\la f_{i_{r}},\g_{\phi(r)}\ra\la f_{i_{r+1}},\g_{\psi(r+1)}\ra...\\
...\la f_{i_{t-1}},\g_{\psi(t-1)}\ra g\end{multline*}
if $(i_{1},i_{2}) = (-1,1)$ and is otherwise zero.\\\\
Now we consider $\Zpgj{1}{t-1}|_{k}(\pi_{\cusp{m}}-1)(\g_{3}-1)...(\g_{t-1}-1)$ in the case $(i_{1},i_{2}) = (-1,1).$ We know by the inductive hypothesis that
\begin{multline}
\Zpgj{1}{t-1}|_{k}(\pi_{\cusp{m}}-1)(\g_{3}-1)...(\g_{t-1}-1)\\
= \left\{ \begin{array}{cc}\Bigl(\displaystyle{\int_{i}^{\pi_{\cusp{m}}i}}\CZ'_{-1, f_{1}}(t)dt\Bigr)\la f_{j_{3}},\g_{2}\ra...\la f_{j_{t-1}},\g_{t-2}\ra g& if (j_{1},j_{2}) = (-1,1)\\
0 & otherwise.\end{array}\right.\label{jhit}
\end{multline}
We also know by the definition of $j_{1},...,j_{t-1}$ that there is a $(\phi,\psi)\in \Shuf{r}{t}$ such that \begin{multline}\la f_{j_{1}},\g_{1}\ra...\la f_{j_{t-1}},\g_{t-1} \ra g =\\ \la f_{i_{1}},\g_{\phi(1)}\ra...\la f_{i_{r}}, \g_{\phi(r)}\ra \la f_{i_{r+1}},\g_{\psi(r+1)}\ra...\la f_{i_{t-1}}, \g_{\psi(t-1)}\ra g\label{jvsi}\end{multline}
for all $\g_{1},...,\g_{t-1}\in \G.$ By the definition of shuffles, this means that $(j_{1},j_{2})$ must be one of $(i_{1},i_{2}),(i_{1}, i_{r+1}), (i_{r+1},i_{1})$ and $(i_{r+1},i_{r+2})$ if $r+1<t-1,$ and one of the first three of these if $r+1 = t-1.$ Since $i_{1}, i_{r+1}, i_{r+2} < 0,$ the case $(j_{1},j_{2}) = (-1,1)$ (and so (\ref{jhit}) is nonzero) can only occur in terms where $(j_{1},j_{2}) = (i_{1},i_{2}),$ and so $\phi(1) = 1,$ $\phi(2) = 2.$ Using this fact, we can cancel the first two terms from each side of (\ref{jvsi}), and the resulting equality lets us rewrite our previous expression as
\begin{multline*}
\CZ'_{-1,1,j_{3},...,j_{t-1};\,g}|_{k}(\pi_{\cusp{m}}-1)(\g_{3}-1)...(\g_{t-1}-1)\\
= \int_{i}^{\pi_{\cusp{m}}i}\CZ'_{-1;\,f_{1}}(t)dt \la f_{i_{3}},\g_{\phi(3)}\ra...\la f_{i_{r-1}}, \g_{\phi(r-1)}\ra\\ \la f_{i_{r}},\g_{\psi(r)}\ra...\la f_{i_{t-1}}, \g_{\psi(t-1)}\ra g.
\end{multline*}
Summing over all the terms of $\mathcal{C}(i_{1},...,i_{t-1};g),$ this gives us a nonzero term in the form given above for each $(\phi,\psi)$ satisfying $\phi(1)=1,$ $\phi(2) = 2$ apart $(\phi,\psi)_{0}$. Subtracting these from \[\Bigl(\Zpg{r+1}{t-1}(z)\int_{i}^{z}\Zpf{1}{r-1}{f_{i_{r}}}(t)dt\Bigr)|_{k}(\pi_{\cusp{m}}-1)(\g_{3}-1)...(\g_{t-1}-1)\] leaves only that term corresponding to $(\phi, \psi)_{0},$ and so if $(i_{1},i_{2})=(-1,1),$
\[\Zpg{1}{t-1}|_{k}(\pi_{\cusp{m}}-1)(\g_{3}-1)...(\g_{t-1}-1) = \int_{i}^{\pi_{\cusp{m}}i}\CZ'_{-1;\,f_{1}}(t)dt \la f_{i_{3}},\g_{3}\ra...\la f_{i_{t-1}},\g_{t-1}\ra g\]

If $(i_{1},i_{2}) \neq (-1,1),$ we know that \[\Bigl(\Zpg{r+1}{t-1}(z)\int_{i}^{z}\Zpf{1}{r-1}{f_{i_{r}}}(t)dt\Bigr)|_{k}(\pi_{\cusp{m}}-1)(\g_{3}-1)...(\g_{t-1}-1) = 0.\]
Further, since $i_{1}, i_{r+1}, i_{r+2} < 0$ and $(i_{1}, i_{2}) \neq (-1,1),$ we know that $(j_{1},j_{2})\neq (-1,1)$ for all terms $\Zpgj{1}{t-1}$ of $\mathcal{C}(i_{1},...,i_{t-1};\,g),$ and so \[\mathcal{C}(i_{1},...,i_{t-1};\,g)|_{k}(\pi_{\cusp{m}}-1)(\g_{3}-1)...(\g_{t-1}-1) = 0.\] Thus $\Zpg{1}{t-1}|_{k}(\pi_{\cusp{m}}-1)(\g_{3}-1)...(\g_{t-1}-1) = 0,$ as desired. The final case, where $i_{1},...,i_{t-1}<0$ follows from the previously stated properties of the series used to define $\Zpg{1}{t-1}$ in these cases.
\end{proof}
%
%
%
%
%
%
\begin{Lemma}
Given $F \in M_{2}^{k}(\G)$ and $\cstuff \in \C$ such that
\[F|_{k}(\g_{1}-1)...(\g_{t-1}-1) = \sum_{(i_{3},...,i_{t-1};\,g)\in J_{t-2,k}} \cstuff \la f_{-1},\g_{1}\ra...\la f_{i_{t-1}},\g_{t-1}\ra g\] for all $\g_{1},...,\g_{t-1}\in\G,$ then we must have $\cstuff = 0$ for all $i_{3},...,i_{t-1}$ and $g.$\label{Fimposs}
\end{Lemma}
\begin{proof}
Given such an $F,$ consider $F- \sum \cstuff \Zpog{3}{t-1}.$ We know that
\[\bigl(F- \sum \cstuff \Zpog{3}{t-1}\bigr)|_{k}(\g_{1}-1)...(\g_{t-1}-1) = 0,\] for all $\g_{1},...,\g_{t-1}\in\G,$ and so
$\bigl(F- \sum \cstuff\Zpog{3}{t}\bigr)|_{k}(\g_{1}-1)...(\g_{t-2}-1)$ will vanish if any further $(\g-1)$ is applied. Furthermore, this expression is clearly holomorphic, and it follows from the definition of $\Zpg{1}{t-1}$ that it has at most polynomial growth at all cusps. Thus we can write
\[\bigl(F- \sum \cstuff\Zpog{3}{t-1}\bigr)|_{k}(\g_{1}-1)...(\g_{t-2}-1) = \sum_{f_{j}\in\CM_{k}}\chi_{j}(\g_{1},...,\g_{t-2})f_{j},\] where $\chi_{j}:\G^{t-1}\to \C.$\\\\
Since $(\g_{i}\d_{i}-1) = (\g_{i}-1)+(\d_{i}-1)+(\g_{i}-1)(\d_{i}-1),$ we see that \begin{multline*}\bigl(F- \sum \cstuff\Zpog{3}{t-1}\bigr)|_{k}(\g_{1}-1)...(\g_{i}\d_{i}-1)...(\g_{t-1}-1)\\ = \bigl(F- \sum \cstuff\Zpog{3}{t-1}\bigr)|_{k}(\g_{1}-1)...(\g_{i}-1)...(\g_{t-1}-1)\\ + \bigl(F- \sum \cstuff\Zpog{3}{t-1}\bigr)|_{k}(\g_{1}-1)...(\d_{i}-1)...(\g_{t-1}-1).\end{multline*}
The linear independance of the $f_{j}$ then means that $\chi_{j}(\g_{1},...,\g_{i}\d_{i},...,\g_{t}) - \chi_{j}(\g_{1},...,\g_{i},...,\g_{t})+ \chi_{j}(\g_{1},...,\d_{i},...,\g_{t}) = 0$ for $1\le i\le t-1$ - in other words, $\chi_{j}$ is a group homomorphism in terms of each argument. Thus repeatedly applying the Eichler-Shimura isomorphism allows us to write
\begin{multline*}\bigl(F- \sum \cstuff\Zpog{3}{t-1}\bigr)|_{k}(\g_{1}-1)...(\g_{t-2}-1)\\ = \sum_{\substack{k_{i}\in \{\pm 1,...,\pm g+m-1\}\\f_{j}\in\CM_{k}}} a_{k_{1},...,k_{t-2},j}\la f_{k_{1}},\g_{1}\ra...\la f_{k_{t-2}}, g_{t-2}\ra f_{j}\end{multline*}
with $a_{k_{1},...,k_{t-2},j}\in \BBC.$ Furthermore, we know that this expression vanishes if any of the $\g_{i} = \pi_{\mathfrak{a}}$ for $\mathfrak{a} \neq \cusp{m},$ and so the only nonzero terms in the above sum must be those with all $f_{k_{i}}$ cuspidal. This in turn means that
\[\bigl(F- \sum \cstuff\Zpog{3}{t-1}\bigr)|_{k}(\pi_{\cusp{m}}-1)(\g_{2}-1)...(\g_{t-2}-1) = 0\]
for all $\g_{2},...,\g_{t-2}.$ Since $F\in M^{t}_{k},$ we know that $F|_{k}(\pi_{\cusp{m}}-1) = 0,$ and so 
\[\bigl(\sum \cstuff\Zpog{3}{t-1}\bigr)|_{k}(\pi_{\cusp{m}}-1)(\g_{2}-1)...(\g_{t-2}-1) = 0.\]
But applying the previous result gives \begin{multline*}\bigl(\sum \cstuff\Zpog{3}{t-1}\bigr)|_{k}(\pi_{\cusp{m}}-1)(\g_{2}-1)...(\g_{t-2}-1) =\\ \int_{i}^{\pi_{\cusp{m}}i}\CZ_{-1;\,f_{1}}(w)dw \sum \cstuff \la f_{i_{3}},\g_{2}\ra...\la f_{i_{t-1}},\g_{t-2}\ra g
=0\end{multline*}
for all $\g_{2},...,\g_{t-2}.$ Thus by the linear independance of the modular symbols, $\cstuff = 0$ for all $i_{3},...,i_{t-1}$ and $g.$
\end{proof}
\subsection{The Main Theorem}
\begin{Theorem}
Let $t\geq 1.$ Then the image of
\[\{\Zng{1}{t-1}:(i_{1},...,i_{t};\,g)\in I_{t,k}\}\]
under the natural projection is a basis for $M^{t}_{k}/M^{t-1}_{k}.$
\end{Theorem}
\begin{proof}
We have already established that the $\Zng{1}{t-1}$ are in $M^{t}_{k}(\G).$ Their linear independence follows simply from their functional equations.\\\\
To show that the $\Zng{1}{t}$ span $M^{t}_{k},$ consider $F \in M^{t}_{k}(\G).$ As in the previous lemma, we see that we can write
\begin{multline*}F|_{k}(\g_{1}-1)...(\g_{t-1}-1) = \sum_{j}\chi_{j}(\g_{1},...,\g_{t-1})g_{j}\\
= \sum_{(i_{1},...,i_{t-1};\,g_{j})\in I_{t,k}}c_{i_{1},...,i_{t-1};\,g_{j}}\la f_{i_{1}}, \g_{1}\ra...\la f_{i_{t-1}},\g_{t-1} \ra g_{j}.
\end{multline*}
We now consider \begin{multline}(F-L)|_{k}(\g_{1}-1)...(\g_{t-1}-1) :=\\
\Bigl(F - \sum_{(i_{1},...,i_{t-1};\,g_{j})\in I_{t,k}}c_{i_{1},...,i_{t-1};\,g_{j}}\Zng{1}{t-1}\Bigr)|_{k}(\g_{1}-1)...(\g_{t-1}-1).\label{zero}\end{multline} 
To establish the theorem, it will suffice to show that this is 0, and thus that $F$ is in the linear span of $\{\Zng{1}{t-1}\}$ modulo $M_{k}^{t-1}(\G).$\\\\
By the functional equation for the $\Zng{1}{t-1},$ we know that the above expression is simply $\phi(\g_{1},...,\g_{t-1})$ for some $\phi \in A_{t,k}.$ In other words, it is a linear combination
\begin{multline}
\phi(\g_{1},...,\g_{t-1})\\
\shoveleft{=\sum_{(i_{1},...,i_{t-1};g)\in I_{t,k}}c_{i_{1},...,i_{t-1};g}\la f_{i_{1}},\g_{1}\ra...\la f_{i_{t-1}}, \g_{t-1}\ra g}\\
\shoveleft{=\sum_{r=1}^{t-1}\:\sum_{\substack{i_{1},...,i_{r-1}\in\{\pm1,...,\pm g\}\\(i_{r+2},...,i_{t-1};g)\in I_{t-r-1,k}}}
c_{i_{1},...,-1,1,...,i_{t-1},g}\la f_{i_{1}},\g_{1}\ra...\la f_{i_{r-1}},\g_{r-1}\ra} \\
\shoveright{\la f_{-1},\g_{r}\ra \la f_{1}, \g_{r+1}\ra\la f_{i_{r+2}}, \g_{r+2}\ra...\la f_{i_{t-1}},\g_{t-1}\ra g}\\
\shoveleft{=\sum_{r=1}^{t-2}\:\sum_{i_{1},...,i_{r-1}\in \{\pm1,...,\pm g\}} a_{i_{1},...,i_{r-1}}\la f_{i_{1}},\g_{1}\ra...\la f_{i_{r-1}},\g_{r-1}\ra \la f_{-1},\g_{r}\ra} \\
\shoveright{\la f_{1}, \g_{r+1}\ra \sum_{(i_{r+2},...,i_{t-1};g)\in I_{t-r-1,k}}b_{i_{r+2},...,i_{t-1};g}\la f_{i_{r+2}}, \g_{r+2}\ra...\la f_{i_{t-1}},\g_{t-1}\ra g}\\
\shoveright{+\sum_{i_{1},...,i_{t-2}\in \{\pm1,...,\pm g\}} a_{i_{1},...,i_{t-2}}\la f_{i_{1}},\g_{1}\ra...\la f_{i_{t-2}},\g_{r-1}\ra \la f_{-1},\g_{r}\ra f_{1} }
\label{write}\\
\end{multline}
where the final term, corresponding to $r = t-1,$ is nonzero only if $k=2.$\\\\
However, we know from (\ref{maintheorem}) that for $(i_{r+2},...,i_{t-1};g) \in I$ there is an $G$ in $M^{t-r}_{k}$ such that 
\begin{multline}G|_{k}(\g_{r+1}-1)...(\g_{t-1}-1) = \la f_{1}, \g_{r+1}\ra\la f_{i_{r+2}},\g_{r+2}\ra...\la f_{i_{t-1}}, \g_{t-1}\ra g\\ + \psi(\g_{r+1},...,\g_{t-1})\label{formform}\end{multline} for some $\psi \in A_{t-r},k.$ Using this we will write
\begin{multline}
\phi(\g_{1},...,\g_{t-1}) = \sum_{r=1}^{t-1}\sum_{i_{1},...,i_{r-1}\in \{\pm1,...,\pm g\}}d_{i_{1},...,i_{r-1}}\la f_{i_{1}},\g_{1}\ra...\la f_{i_{r-1}}, \g_{r-1}\ra\\
\la f_{-1}, \g_{r}\ra\Bigl(F_{i_{1},...,i_{r-1}}|_{k}(\g_{r+1}-1)...(\g_{t-1}-1)\Bigr).\label{rewrite}\end{multline}
We do this by starting with $r = 1,$ and substituting (\ref{formform}) to rewrite the terms in (\ref{write}) corresponding to each given $r$ as $r$ increases - we can treat the terms $\la f_{i_{1}},\g_{1}\ra...\la f_{i_{r-1}},\g_{r-1}\ra \la f_{-1},\g_{r}\ra\phi(\g_{r+1},...,\g_{t-1})$ that are left over along with the terms of (\ref{write}) corresponding to higher $r$, as it will be of the same form as them by the definition of $\mathcal{A}_{t-r,k}.$ For a given $(i_{1},...,i_{r-1})$ we can then write $F_{i_{1},...,i_{r-1}}$ for the sum over all $i_{r+1},..., i_{-1},g$ of the corresponding $G.$ The final application of (\ref{formform}) will yield no left over $\phi,$ since the forms being substituted are simply of the form $G(z) = g(z)\int_{i}^{z}f_{1}(t)dt$ for $k>2$ or simpy $G(z) = f_{1}(z)$ if $k=2.$\\\\
Recall that $\phi(\g_{1},...,\g_{t-1}) = (F-L)|_{k}(\g_{1}-1)...(\g_{t-2}-1)(\g_{t-1}-1).$ If $k=2$ we take all terms except the last one over to the left hand side of (\ref{rewrite}) to give
\begin{multline*}
\Bigl[\bigl((F-L)|_{k}(\g_{1}-1)...(\g_{t-2}-1)\bigr)\Bigl.\\ 
-\sum_{r=1}^{t-2}\sum_{i_{1},...,i_{r-1}\in \{\pm1,...,\pm g\}}d_{i_{1},...,i_{r-1}}\la f_{i_{1}},\g_{1}\ra...\la f_{i_{r-1}}, \g_{r-1}\ra\\
\Bigr.\la f_{-1},\g_{r}\ra F_{i_{1},...,i_{r-1}}|_{k}(\g_{r+1}-1)...(\g_{t-2}-1)\Bigr]|_{k}(\g_{t-1}-1)\\
= \Bigl(\sum_{(i_{1},...,i_{t-2})}d_{i_{1},...,i_{t-2}}\la f_{i_{1}},\g_{1}\ra...\la f_{i_{t-2}},\g_{2}\ra\Bigr)\cdot \la f_{-1}, \g_{t-1}\ra f_{1},\end{multline*}
and because for fixed $\g_{1},...,\g_{t-2}$ the right hand side is a multiple of $\la f_{-1},\g_{t-1}\ra f_{1}$ and the terms in the square brackets on the LHS are all in $M^{2}_{k}(\G),$  (5.2) of \cite{DO} tells us that its coefficient must be zero. Since this holds for any $\g_{1},...,\g_{t-2}$ and the products of modular forms are linearly independent, this means that each $d_{i_{1},...,i_{t-2}}$ is zero. For $k>2,$ this term is zero automatically.\\\\
We then apply the same process inductively, proving that terms are zero in order of decreasing $r.$ At each stage, having previously shown that any terms with $r>n$ are zero, we fix $\g_{1},...,\g_{n-1}$ and take the terms with $r<n$  over to the left hand side. The remaining terms - those with $r=n$ - can also be written in the form \begin{multline*}\sum_{i_{1},...,i_{n-1}\in \{\pm1,...,\pm g\}} a'_{i_{1},...,i_{n-1}}\la f_{i_{1}},\g_{1}\ra...\la f_{i_{n-1}},\g_{n-1}\ra\cdot\\
\sum_{(i_{n+2},...,i_{t-1};g)\in I_{t-n-1,k}}b'_{i_{n+2},...,i_{t-1};g}\la f_{-1},\g_{n}\ra\la f_{1}, \g_{n+1}\ra\la f_{i_{n+2}}, \g_{n+2}\ra...\la f_{i_{t-1}},\g_{t-1}\ra g\end{multline*}
where the $a'$ and $b'$ coming from the $a$ and $b$ in (\ref{write}) and from the additional $\phi$ terms added in during the rewriting process above. In other words, we have
\begin{multline*}
\Bigl[\bigl((F-L)|_{k}(\g_{1}-1)...(\g_{n-1}-1)\bigr)\Bigl.\\ 
-\sum_{r=1}^{n-2}\sum_{i_{1},...,i_{r-1}\in \{\pm1,...,\pm g\}}d_{i_{1},...,i_{r-1}}\la f_{i_{1}},\g_{1}\ra...\la f_{i_{r-1}}, \g_{r-1}\ra\\
\Bigr.\la f_{-1},\g_{r}\ra F_{i_{1},...,i_{r-1}}|_{k}(\g_{r+1}-1)...(\g_{n-1}-1)\Bigr]|_{k}(\g_{n}-1)...(\g_{t-1}-1)\\
\shoveleft{=\sum_{i_{1},...,i_{n-1}\in \{\pm1,...,\pm g\}} a'_{i_{1},...,i_{n-1}}\la f_{i_{1}},\g_{1}\ra...\la f_{i_{n-1}},\g_{n-1}\ra\cdot } \\
\shoveright{ \sum_{(i_{n+2},...,i_{t-1};g)\in I_{t-n-1,k}}b'_{i_{n+2},...,i_{t-1};g}\la f_{-1},\g_{n}\ra\la f_{1}, \g_{n+1}\ra\la f_{i_{n+2}}, \g_{n+2}\ra...\la f_{i_{t-1}},\g_{t-1}\ra g.}\end{multline*}
Thus we can apply the previous lemma, which, on varying $\g_{1},...,\g_{n-1},$ tells us that each $b'_{i_{n+2},...,i_{t-1};g}$ is zero. Doing this for every $r,$ we conclude that
\[(F-L)|_{k}(\g_{1}-1)...(\g_{t-1}-1) = 0.\]
\end{proof}
%
%
%
%
%
%
%
%
%
%
%
\section{Cohomology}
\subsection{The Spaces $\ZM{t}$ and $\ZS{t}$}
\begin{Definition}\begin{align*}
\ZM{t} = \bigoplus_{r=0}^{t-1}\bigl\la \CZ'_{i_{1},...,i_{r-1};\,f}|i_{1},...,i_{r}\in \{\pm1,...,\pm g\},\: f\in M_{k}(\G) \bigr\ra\\
\ZS{t} = \bigoplus_{r=0}^{t-1}\bigl\la \CZ'_{i_{1},...,i_{r-1};\,f}|(i_{1},...,i_{r})\in \{\pm1,...,\pm g\},\: f\in S_{k}(\G) \bigr\ra
\end{align*}
\end{Definition}
For $t=1,$ these are simply $M_{k}(\G)$ and $S_{k}(\G).$ In general, we have
\begin{eqnarray*}
\dim_{\BBC}\ZM{t}& =& \bigl((2g)^{t}+...+2g+1\bigr)\dim_{\BBC}M_{k}(\G)\\
\dim_{\BBC}\ZS{t}& =& \bigl((2g)^{t}+...+2g+1\bigr)\dim_{\BBC}S_{k}(\G).
\end{eqnarray*}
\begin{Lemma}$\ZM{t}$ is precisely the space of functions $f: \mathbb{H} \to \mathbb{C}$ satisfying\\
\begin{itemize}
\item f is holomorphic.
\item For each $r \geq 0$ and $\g_{1},...,\g_{r}\in\G,$ $f|_{k}(\g_{1}-1)...(\g_{r}-1)$ has at most polynomial growth at the cusps.
\item $f|_{k}(\pi_{\cusp{i}} - 1) = 0$ for all $i \neq m.$
\item $f|_{k}(\g_{1} - 1)...(\g_{t}-1) = 0$ for all $\g \in \G.$
\end{itemize}
$\ZS{t}$ is the space of functions satisfying these conditions as well as
\begin{itemize}
\item For each $r \geq 0$ and $\g_{1},...,\g_{r}\in\G,$ $f|_{k}(\g_{1}-1)...(\g_{r}-1)$ decays exponentially at the cusps.
\end{itemize}
\end{Lemma}
\begin{proof}
In both cases, this is true by definition for $t = 1.$ Now suppose that the proposition holds for $\ZM{n}$ for all $n<t.$ If we write $D^{(t)}$ for the space satisfying the first four conditions above, we can see that $\ZM{t} \hookrightarrow D^{(t)}$ - we have already checked all but the second condition, and that follows directly from the properties of the series used to construct $\CZ'_{i_{1},...,i_{r-1};\,f}.$\\\\
Now, for $f \in D^{(t)},$ we can see that $f|_{k}(\g_{1}-1)...(\g_{t-1}-1) \in M_{k}(\G),$ and so we can write
\[f|_{k}(\g_{1}-1)...(\g_{t-1}-1) = \sum_{g_{i}\in\CM_{k}}\chi_{i}(\g_{1},...,\g_{t-1})g_{i}\] for some functions $\chi_{i}:\G^{t-1}\to \C.$ Since $(\g\d-1) = (\g-1)(\d-1) +(\g-1) + (\d-1),$ we see that the $\chi_{i}$ must be homomorphisms on each term. Thus
\[f|_{k}(\g_{1}-1)...(\g_{t-1}-1) = \sum_{i}\Bigl(\sum_{j_{1},...,j_{t}}c_{j_{1},...,j_{t-1}}^{i}\la f_{j_{1}},\g_{1} \ra \la f_{j_{t-1}},\g_{t-1} \ra\Bigr)g_{i}\label{modsyms}\]
with $c^{i}_{j_{1},...,j_{t-1}}\in \C$ and $f_{j_{k}}\in\CM_{2}$ Furthermore, since $(\g-1)(\pi_{\cusp{i}}-1) = (\g\pi_{\cusp{i}}\g^{-1}-1)\g - (\pi_{\cusp{i}} - 1),$ the left hand side of (\ref{modsyms}) vanishes if any of the $\g$ is $\pi_{\cusp{i}}$ for $i < m,$ and so each of the $f_{j}$ must be cuspidal.\\\\
We now define the map $\varphi: D^{(t)} \to M_{k}(\G)^{(2g)^{t-1}}$ by setting \[\varphi(f) = (f|_{k}(\g_{1}-1)...(\g_{t-1}-1))_{\g_{1},...,\g_{t-1}\in\G}.\]
Since $\mathbf{Z}_{M}^{(t)}\subseteq D^{(t)},$ we know that $\varphi$ is surjective. Its kernel is precisely the set of maps in $D^{(t)}$ which satisfy $f|_{k}(\g_{1}-1)...(\g_{t-1}-1) = 0$ for all $\g_{1},...,\g_{t-1}\in \G,$ which is $D^{(t-1)}$ by definition and thus $\mathbf{Z}_{M,k}^{(t-1)}$ by induction. Thus we have a short exact sequence
\[0\to \mathbf{Z}_{M,k}^{(t-1)} \to D^{(t)} \stackrel{\varphi}{\to} M_{k}(\G)^{(2g)^{t-1}} \to 0,\]
and a comparison of dimensions then establishes that $D^{(t)} = \ZM{t},$ as desired.\\\\
An identical argument works in the case of $\ZS{t}.$ 
\end{proof}
\begin{Corollary}
If $f \in \ZM{t}$ (resp. $\ZS{t}$) then $f|_{k}(\g-1) \in \ZM{t-1}$ (resp. $\ZS{t-1}$).\end{Corollary}
\begin{proof}
For such $f,$ $f|_{k}(\g-1)(\pi_{\cusp{i}}-1) = f|_{k}(\g\pi_{\cusp{i}}-1)\g^{-1} - (\pi_{\cusp{i}}-1) = 0.$ The other conditions given above for $f|_{2}(\g-1)$ to be in $\ZM{t-1}$ or $\ZS{t-1}$ follow immediately from the conditions for $f$ to be in $\ZM{t}$ or $\ZS{t}.$
\end{proof}
\subsection{Construction of Cohomology Groups}
\begin{Theorem}
 For each $t>1,$ there exist subspaces $\Z{t}$ and $\B{t}$ of $C^{1}(\G, P_{k-2})$ such that the homomorphism
\[\psi:\ZM{t}\to C^{1}(\G,P_{k-2})\] which takes $F\in \ZM{t}$ to the map 
$\psi_{F}:\g\mapsto \int_{i}^{\g^{-1}i}F(z)(z-X)^{k-2}dz$ induces an injection
\[\psi': \frac{\overline{\ZS{t}}\bigoplus \ZM{t}}{\overline{\ZS{t-1}}\bigoplus \ZM{t-1}}\hookrightarrow \Hom{t}:= \frac{\Z{t}}{\B{t}}.\]
\end{Theorem}
\begin{proof}
 Set
\[\Z{0} = Z^{1}(\G, P_{k-2}),\] \[\B{0} = B^{1}(\G, P_{k-2})\] and \[\Z{n} = \B{n} = \{0\}\] for $n < 0.$\\\\
We define $\alpha:C^{1}(\G,P_{k-2})\to C^{1}(\G,C^{1}(\G, P_{k-2}))$ by setting \[\alpha(\psi)(\g)(\d)=d\psi(\g, \d)|_{2-k}\g^{-1}.\]
We now work inductively.\\\\
Suppose that for each $n<t$ there exist $\Z{n}$ and $\B{n}$ in $C^{1}(\G,P_{k-2})$ such that
\begin{enumerate}
 \item  $\Z{n-1}\subset \B{n},$ \label{hyp1}
\item for any $\psi \in \Z{n}$ (resp. $\B{n}$) and $\g \in \G,$ $\alpha(\psi)(\g)\in\Z{n-1}$ (resp $\B{n-1}$) and \label{hyp2}
\item the homomorphism
\[\psi:\ZM{n}\to C^{1}(\G,P_{k-2})\] which takes $F\in \ZM{n}$ to the map 
$\psi_{F}:\g\mapsto \int_{i}^{\g^{-1}i}F(z)(z-X)^{k-2}dz$ induces an injection
\[\psi': \frac{\overline{\ZS{n}}\bigoplus \ZM{n}}{\overline{\ZS{n-1}}\bigoplus \ZM{n-1}}\hookrightarrow \Hom{n}:= \frac{\Z{n}}{\B{n}}.\] \label{hyp3}
\end{enumerate}
For $t=2$ this is done in \cite{DO}, the final statement following from the classical Eichler-Shimura theorem by exact analogy with Theorem 7.1 of that paper.\\\\
Now, writing $\overline{\cdot}$ for the reduction map \[\pi:C^{1}(\G,P_{k-2}) \to \frac{C^{1}(\G,P_{k-2})}{\Z{t-2}},\]
we can define
\[\alpha_{t}:C^{1}(\G,P_{k-2}) \to C^{1}\Bigl(\G,\frac{C^{1}(\G,P_{k-2})}{\Z{t-2}}\Bigr)\] by setting $\alpha_{t}(\psi)(\g) = \overline{\alpha(\psi)(\g)}.$\\\\
If we use the trivial action of $\G$ on $C^{1}(\G, P_{k-2})$ then applying \ref{hyp1}, we see that
\[H^{1}_{\Par}\Bigl(\G, \frac{\Z{t-1}}{\Z{t-2}}\Bigr) = Z^{1}_{\Par}\Bigl(\G, \frac{\Z{t-1}}{\Z{t-2}}\Bigr) \subset C^{1}\Bigl(\G,\frac{C^{1}(\G,P_{k-2})}{\Z{t-2}}\Bigr)\]
and
\[H^{1}_{\Par}\Bigl(\G, \frac{\B{t-1}}{\Z{t-2}}\Bigr) = Z^{1}_{\Par}\Bigl(\G, \frac{\B{t-1}}{\Z{t-2}}\Bigr) \subset C^{1}\Bigl(\G,\frac{C^{1}(\G,P_{k-2})}{\Z{t-2}}\Bigr),\]
and so we can define
\begin{eqnarray*}
\Z{t}&=&\alpha_{t}^{-1}H^{1}_{\Par}\Bigl(\G, \frac{\Z{t-1}}{\Z{t-2}}\Bigr),\\
\B{t}&=&\alpha_{t}^{-1}H^{1}_{\Par}\Bigl(\G, \frac{\B{t-1}}{\Z{t-2}}\Bigr),\\
\end{eqnarray*}
and set
\[\Hom{t} = \frac{\Z{t}}{\B{t}}.\]
Now, taking $\psi\in \Z{t-1},$ we know by \ref{hyp2} that $\alpha(\psi)(\g)\in\Z{t-2}.$ But this means that $\overline{\alpha(\psi)(\g)} = 0$ for all $\g,$ and so $\alpha_{t}(\psi) = 0.$ Thus $\psi \in \B{t},$ and so $\Z{t-1}\subset\B{t}$ - condition $\ref{hyp1}$ is satisfied for $n= t.$ Furthermore, since for $\psi \in \Z{t}$ we know by definition that
\[\alpha_{t}(\psi)(\g) = \overline{\alpha(\psi)(\g)} \in \frac{\Z{t-1}}{\Z{t-2}}\]
and so we must have $\alpha(\psi)(\g) \in \Z{t-1}.$ The same is true, mutatis mutandis, for $\B{t},$ and so \ref{hyp2} is also satisfied for $n=t.$\\\\
A little more work is required for \ref{hyp3}. Given $F \in \overline{\ZS{t}}\bigoplus \ZM{t},$ we define
\[\psi_{F}(\g) = \int_{i}^{\g^{-1}i}F(z)(z-X)^{k-2}dz \in C^{1}(\G, P_{k-2}).\]
A simple calculation shows that
\[\alpha(\psi_{F})(\g)(\d) = d\psi_{f}(\g)(\d)|_{2-k}(\g^{-1}) = \int_{i}^{\d^{-1}i}F|_{k}(\g^{-1}-1)(z)(z-X)^{k-2}dz,\]
and for a fixed $\g,$ this is in $\Z{t-1}$ by the inductive hypothesis. Furthermore,
\begin{multline*}\alpha(\psi_{F})(\g_{1}\g_{2})(\d) = \alpha(\psi_{F})(\g_{1})(\d) + \alpha(\psi_{F})(\g_{2})(\d)\\+ \int_{i}^{\d^{-1}i}F|_{k}(\g_{1}^{-1}-1)(\g_{2}^{-1}-1)(z)(z-X)^{k-2}dz
 \end{multline*}
Applying the projection $\pi$ to each side of this, we find that $\alpha_{t}(\psi_{F})(\g_{1}\g_{2})= \alpha_{t}(\psi_{F})(\g_{1}) + \alpha_{t}(\psi_{F})(\g_{2}),$ and so \[\alpha_{t}(\psi_{F}) \in H_{\Par}^{1}\Bigl(\G,\frac{\Z{t-1}}{\Z{t-2}}\Bigr).\] Moreover, for $F\in  \overline{\ZS{t-1}}\bigoplus \ZM{t-1},$ $\alpha(\psi_{F})(\g)\in \Z{t-2}$ for all $\g,$ and so \[\alpha_{t}(\psi_{F}) = 0 \in H_{\Par}^{1}\Bigl(\G,\frac{\B{t-1}}{\Z{t-2}}\Bigr).\]
Thus $\psi_{F}\in \Z{t}$ for $F\in  \overline{\ZS{t}}\bigoplus \ZM{t}$ and $\psi_{F}\in \B{t}$ for $F\in  \overline{\ZS{t-1}}\bigoplus \ZM{t-1}.$\\\\
We now consider the map 
\[\psi': \frac{\overline{\ZS{t}}\bigoplus \ZM{t}}{ \overline{\ZS{t-1}}\bigoplus \ZM{t-1}}\to \Hom{t}\]
induced by $\psi.$ To see that it is well defined, we note that if $[F] = [G]$ (where we are using $[\cdot]$ to represent equivalence classes in the quotient) then $F-G \in  \overline{\ZS{t-1}}\bigoplus \ZM{t-1},$ and so $\psi_{F}-\psi_{G} = \psi_{F-G} \in \Z{t-1}$ by the inductive hypothesis. Since $\Z{t-1} \subset \B{t},$ this means that $\psi_{F} = \psi_{G}$ in $\Hom{t}.$\\\\
For injectivity, take $F \in  \overline{\ZS{t}}\bigoplus \ZM{t}$ such that $\psi'([F]) = 0.$ Then $\psi_{F} \in \B{t},$ and so for all $\g,$ 
\[\alpha(\psi_{F})(\g,\d) = \int_{i}^{\d^{-1}i}F|_{k}(\g^{-1}-1)(z)(z-X)^{k-2}dz \in \B{t-1}\]
as a function of $\d.$\\\\
But the above expression shows that \[\alpha(\psi(F))(\g)(\d) = \psi_{F|_{k}(\g^{-1}-1)}(\d),\] and since $F|_{k}(\g^{-1}-1) \in  \overline{\ZS{t-1}}\bigoplus \ZM{t-1},$ and \ref{hyp3} tells us that $\psi'$ is injective in the case $n=t-1,$ we see that if $\psi_{F|_{k}(\g^{-1}-1)}(\d)$ is a coboundary then $F|_{k}(\g^{-1}-1)$ must in fact be in $ \overline{\ZS{t-2}}\bigoplus \ZM{t-2}$ for all $\g.$  Thus \[F \in  \overline{\ZS{t-1}}\bigoplus \ZM{t-1}\] and so $[F]=[0].$
\end{proof}
\subsection{An Eichler-Shimura map for higher order modular forms}
\begin{Lemma}
The map $\alpha_{t}$ induces an injection
\[\Hom{t}\hookrightarrow \frac{H^{1}_{\Par}(\G, \frac{\Z{t-1}}{\Z{t-2}})}{H^{1}_{\Par}(\G, \frac{\B{t-1}}{\Z{t-2}})}\]
\end{Lemma}
\begin{proof}
Both the well definedness and injectivity of the derived map follow from the definition of $\B{t}.$
\end{proof}

\begin{Proposition}
 \[\frac{H^{1}_{\Par}(\G, \frac{\Z{t-1}}{\Z{t-2}})}{H^{1}_{\Par}(\G, \frac{\B{t-1}}{\Z{t-2}})}\cong H^{1}_{\Par}(\G, \Hom{t-1}).\]
\end{Proposition}
\begin{proof}
 This proof exactly follows the model of Theorem 7.1 of \cite{DO}. From the short exact sequence
\[0 \to \B{t-1}\hookrightarrow \Z{t-1} \to \Hom{t-1} \to 0\]
we derive the long exact sequence
\begin{multline*}
 H^{1}_{\Par}(\G, \frac{\B{t-1}}{\Z{t-2}})\to H^{1}_{\Par}(\G,\frac{\Z{t-1}}{\Z{t-2}})\\
\to H^{1}_{\Par}(\Hom{t-1})\to H^{2}_{\Par}(\G, \frac{\B{t-1}}{\Z{t-2}})\\
\to H^{2}_{\Par}(\G,\frac{\Z{t-1}}{\Z{t-2}})\to H^{2}_{\Par}(\Hom{t-1})\\
\to H^{3}(\G, \frac{\B{t-1}}{\Z{t-2}}) \to \dots.
\end{multline*}
As in \cite{DO}, we know that $H^{j}(\G,M) = 0$ for every $j\geq2$ and for every $\C-$vector space $M$, and that $H^{2}_{par}(\G, M)\cong M/M_{1},$ where $M_{1}$ is the subspace of $M$ generated by the elements $M.(\g-1)$ for $\g\in\G.$ In this case, $M_{1}$ is trivial since we are using the trivial action of $\G$ on the relevant spaces, and so we get the exact sequence
\begin{multline*}
 0\to i^{*}\Bigl(H^{1}_{\Par}(\G, \frac{\B{t-1}}{\Z{t-2}})\Bigr)\to H^{1}_{\Par}(\G,\frac{\Z{t-1}}{\Z{t-2}})\\
\to H^{1}_{\Par}(\G, \Hom{t-1}) \to \B{t-1}\\ \to \Z{t-1} \to \Hom{t-1}\to 0
\end{multline*}
where $i^{*}$ is the map derived from the injective term of the original short exact sequence.\\\\
When $\G$ acts trivially on $M,$ as is the case for all of our coefficient modules, $B^{1}_{\Par}(\G, M)$ is trivial by definition. This means that that $i^{*}$ is injective and so the first term can be replaced with $H^{1}_{\Par}(\G, \frac{\B{t-1}}{\Z{t-2}}).$ Counting dimensions (and bearing in mind that the last three terms here are the terms of the original short exact sequence), this tells us that
\begin{multline*}0\to H^{1}_{\Par}\Bigl(\G, \frac{\B{t-1}}{\Z{t-2}}\Bigr)\to H^{1}_{\Par}\Bigl(\G,\frac{\Z{t-1}}{\Z{t-2}}\Bigr)\\
\to H^{1}_{\Par}(\G, \Hom{t-1})\to0\end{multline*} is exact, and the result follows.
\end{proof}

Now, if we suppose the inductive hypothesis that $\Hom{t-1}\cong \ZFrac{t-1}{t-2},$ we can compose our maps to get an injection
\begin{eqnarray*}
 \ZFrac{t}{t-1}&\hookrightarrow&\Hom{t}\\
&\hookrightarrow&\frac{H^{1}_{\Par}(\G, \frac{\Z{t-1}}{\Z{t-2}})}{H^{1}_{\Par}(\G, \frac{\B{t-1}}{\Z{t-2}})}\\
&\hookrightarrow&H^{1}_{\Par}(\G, \Hom{t-1})\\
&\tilde{\rightarrow}&\bigoplus_{i=1}^{2g}\bigl(\ZFrac{t-1}{t-2}\bigr)\\
\end{eqnarray*}
Viewing them as the linear spans of the $\CZ_{i_{1},...,i_{t-1};\,f},$ we know that the dimensions of $\ZM{t}$ and $\ZS{t}$ to be $(2g)^{t-1}\dim M_{k}(\G)$ and $(2g)^{t-1}\dim S_{k}(\G)$ respectively. Thus, by comparing dimensions in the above chain of injections, we see that they must also be surjections. Thus we have established
\begin{Theorem}
$\psi': \ZFrac{t}{t-1} \rightarrow \Hom{t}$ is an isomorphism.
\end{Theorem}
We know that $\bar{S}_{k}^{t}(\G)\oplus M_{k}^{t}(\G) \subseteq \overline{\ZS{t}}\oplus\ZM{t},$
and we can see from the definitions that
\[\bigl(\overline{\ZS{t-1}}\oplus\ZM{t-1}\bigr)\cap \bigl(\bar{S}_{k}^{t}(\G)\oplus M_{k}^{t}(\G)\bigr) \cong \bar{S}_{k}^{t-1}(\G)\oplus M_{k}^{t-1}(\G).\]
This implies
\begin{Corollary}The map $\psi$ also induces an injection 
\[\psi'': \frac{\bar{S}_{k}^{t}(\G)}{M_{k}^{t}(\G)}\oplus\frac{\bar{S}_{k}^{t}(\G)}{ M_{k}^{t}(\G)}\hookrightarrow \Hom{t}\]
\end{Corollary}

\end{document}